\documentclass{article}

\usepackage{amsmath}
\usepackage{url}
\def\qq{\qquad}

\def\vf{\varphi}

\begin{document}

\title{A Short Proof of Stirling's Formula}
\author{Hongwei Lou}
\date{}
\maketitle

\begin{abstract}
By changing variables in a suitable way and using dominated convergence methods, this note gives a short proof of Stirling's formula and its refinement.
\end{abstract}
% Mathematical Subject Classification: 26A12

\section{Introduction.}

This note provides a short proof of the well-known Stirling's formula
\begin{equation}\label{E0}
\Gamma (s+1)=\Big({s\over e}\Big)^s\sqrt{2\pi s} \, (1+o(1)), \qquad \mbox{as}\, s\to +\infty,
\end{equation}
where $\Gamma$ is Euler's Gamma function
\begin{equation}\label{E0a}
\Gamma (s): =\int^{+\infty}_0 t^{s-1}e^{-t}\, dt, \qquad \mbox{for}\, s>0.
\end{equation}
There are many alternative proofs of Stirling's formula in the literature. For instance, one can find proofs given by P. Diaconis, D. Freedman \cite{D-F},  J. M. Partin  \cite{P} and   M. Pinsky \cite{Pin}. The proof presented below has some common elements with the ones in \cite{D-F} and \cite{P}, but appears much simpler in that it also gives a short and clear way to refine Stirling's formula (cf. (\ref{E15a}) and (\ref{E17})).

\section{Proof of Stirling's formula.}

Let $s>0$. Changing variables ($t=s(x+1)$) leads to
\begin{equation}\label{E1}
\Gamma (s+1)=\int^{+\infty}_0 t^se^{-t}\, dt
={s^{s+1}\over e^s}\int^{+\infty}_{-1}e^{-sg(x)}\, dx,
\end{equation}
where
\begin{equation}\label{E2}
g(x):=x- \ln (1+x), \qquad \mbox{for}\, x\in (-1,+\infty).
\end{equation}
One can easily see that $g(x)$ is strictly convex in $(-1,+\infty)$, admitting its minimum $0$ at $x=0$. Moreover, $\displaystyle\lim_{x\to -1^+}g(x)=\lim_{x\to +\infty}g(x)=+\infty$.

Let $v(x)= \sqrt{g(x)}\,{\rm sgn}\, x$ and $x=f(v)$ be the inverse function of $v(x)$.
By a change of variables ($\displaystyle x=f({u\over \sqrt s})$), one has
\begin{eqnarray}\label{E5}
\nonumber \Gamma(s+1)&=&  \Big({s\over e}\Big)^s\sqrt s\int^{+\infty}_{-\infty} 2\Big({{u\over \sqrt s}\over f({u\over \sqrt s})}+{u\over \sqrt s}\Big) e^{-u^2} \, du\\
&=&  \Big({s\over e}\Big)^s\sqrt s\int^{+\infty}_02\Big(y\big({u\over \sqrt s}\big)+y\big(-{u\over \sqrt s}\big)\Big) e^{-u^2} \, du,
\end{eqnarray}
where
\begin{equation}\label{E6a}
y(v)=\left\{\begin{array}{ll}\displaystyle
 {v\over f(v)}, & v\ne 0,\vspace{4mm}\\
\displaystyle {\sqrt 2\over 2}, & v=0.\end{array}\right.
\end{equation}
Since $f(v)$ is continuous and
\begin{equation}\label{E3}
\lim_{v\to 0} y(v)=\lim_{v\to 0} {v\over f(v)}=\lim_{x\to 0} {\sqrt{x-\ln (1+x)}\over |x|}={\sqrt 2\over 2},
\end{equation} it follows that $y(v)$ is continuous.
Note that $f(v)$ is strictly increasing on  $(-\infty,+\infty)$ with $f(0)=0$, the equations (\ref{E6a}) and (\ref{E3}) imply that
\begin{equation}\label{E4}
0<y(v)\leq C(|v| +1), \qquad\mbox{for all}\,  v\in (-\infty,+\infty)
\end{equation}
for some constant $C>0$.
Thus, one deduces from (\ref{E4}) and the continuity of $y(v)$ at $v=0$ that
\begin{equation}\label{E9}
\lim_{s\to +\infty} \int^{+\infty}_02\Big(y\big({u\over \sqrt s}\big)+y\big(-{u\over \sqrt s}\big)\Big) e^{-u^2} \, du=\int^{+\infty}_04y(0) e^{-u^2}\, du=\sqrt{2\pi},
\end{equation}
which yields Stirling's formula (\ref{E0}).

More precisely, we can prove (\ref{E9}) by virtue of the following estimate:
\begin{eqnarray}\label{E12a}
\nonumber &&  \Big|\int^{+\infty}_0\Big(y\big({u\over \sqrt s}\big)+y\big(-{u\over \sqrt s}\big)-2y(0)\Big) e^{-u^2} \, du\Big|\\
\nonumber &\leq &  \int^{+\infty}_{\sqrt [4] s}\Big(2C\big({u\over \sqrt s}+1\big)+\sqrt 2\Big) e^{-u^2} \, du+\int^{\sqrt [4] s}_02\omega\big({1\over \sqrt [4] s}\big) e^{-u^2} \, du\\
&\leq & \int^{+\infty}_{\sqrt [4] s} (4C+2)u e^{-u^2} \, du+\int^{+\infty}_02\omega\big({1\over \sqrt [4] s}\big) e^{-u^2} \, du, \quad\mbox{for all}\, s\geq 1,
%\nonumber &=& (2C+1)e^{-\sqrt s}+ \omega\big({1\over \sqrt [4] s}\big)\sqrt \pi, \q \mbox{for all}\, s\geq 1
\end{eqnarray}
where
\begin{equation}\label{E9b}
\displaystyle \omega(r):=\max_{|v|\leq r} |y(v)-y(0)|, \qquad  \mbox{for}\, r>0,
\end{equation}
satisfying $\displaystyle \lim_{r\to 0^+}\omega(r)=0$.

\section{Refinement of Stirling's formula.}

To refine Stirling's formula, we mention that
\begin{equation}\label{E10a}
G(x):=\left\{\begin{array}{ll} {g(x)\over x^2}, & x\ne 0, \\
{1\over 2}, & x=0\end{array}\right.%\qquad x\in (-1,+\infty)
\end{equation}
is a positive smooth function in $(-1,+\infty)$ and $y(v)$ is the unique solution of
\begin{equation}\label{E10}
y(v)=\sqrt{G\Big({v\over y(v)}\Big)},  \qquad  \mbox{for}\, v\in (-\infty,+\infty).
\end{equation}
Then the smoothness of $y(v)$ follows easily from the smoothness of $\sqrt G$ and
\begin{equation}\label{E13a}
\displaystyle \Big[{\partial \over \partial y} \Big(y-\sqrt{G\big({v\over y}\big)}\Big)\Big]\Big|_{y=y(v)}\\%={1\over 2(1+{v\over y(v)})G\big({v\over y(v)}\big)}
={1\over 2(1+f(v))y^2(v)}>0.
\end{equation}
Denote by $\displaystyle \sum^\infty_{k=0} a_kv^k$ the MacLaurin's series of $y(v)$. Thus, we have
\begin{equation}\label{E13}
y(v)+y(-v)=\sum^n_{k=0} 2a_{2k}v^{2k}+o(v^{2n}),\qquad  \mbox{for}\,  v\to 0.
\end{equation}
Thus, similar to (\ref{E4}), there exist constants $C_n>0$  (for $n=1,2,\ldots$) such that
\begin{eqnarray}\label{E14}
  && \Big| y(v)+y(-v)-\sum^n_{k=0} 2a_{2k}v^{2k}\Big|\leq C_n v^{2n}, \quad\mbox{for all}\, v\in (-\infty,\infty).
\end{eqnarray}
While analogous to (\ref{E12a}), we have that for any $s\geq 1$,
\begin{eqnarray}\label{E15b}
\nonumber && \Big|\int^{+\infty}_0\Big(y\big({u\over \sqrt s}\big)+y\big(-{u\over \sqrt s}\big)-\sum^n_{k=0} 2a_{2k}{u^{2k}\over s^k}\Big) e^{-u^2} \, du\Big|\\
&\leq & {1\over s^n}\Big[ C_n\int^{+\infty}_{\sqrt [4]s} u^{2n}e^{-u^2}\, du+ \omega_n\big({1\over \sqrt [4]s}\big)\int^{+\infty}_0 u^{2n}e^{-u^2}\, du\Big]
\end{eqnarray}
with
$$
\omega_n(r):=\sup_{0<v\leq r}{1\over v^{2n}}\Big|y(v)+y(-v)-\sum^n_{k=0}a_{2k}v^{2k}\Big|, \qquad \mbox{where}\, r>0.
$$
As a consequence of (\ref{E15b}) and  $\displaystyle \lim_{r\to 0^+}\omega_n(r)=0$,
\begin{eqnarray}\label{E15a}
\nonumber \Gamma(s+1)&=&  \Big({s\over e}\Big)^s\sqrt s\Big[\int^{+\infty}_0  \sum^n_{k=0} 4a_{2k}{u^{2k}\over s^k}e^{-u^2}\, du+o\big({1\over s^n}\big)\Big]\\
&=& \Big({s\over e}\Big)^s\sqrt s\Big[ \sum^n_{k=0} 2a_{2k} \Gamma(k+{1\over 2})  {1\over s^k}  +o\big({1\over s^n}\big)\Big],\quad \mbox{as}\, s\to +\infty.
\end{eqnarray}
For example, for $n=5$, a direct calculation shows that
\begin{eqnarray}
\nonumber\label{E16}
&& a_0={\sqrt 2\over 2}, %\q a_1=-{1\over 3},
\; a_2={\sqrt 2\over 12}, %\q a_3=-{4\over 135},
\; a_4={\sqrt 2\over 432}, %\q a_5={4\over 2835},
\; a_6=-{139\sqrt 2\over 194400}, %\q a_7={8\over 25515},
\\
&& \; a_8=-{571\sqrt 2\over 32659200}, %\q a_9=-{1124\over 37889775},
\; a_{10}={163879\sqrt 2\over 12345177600}.
\end{eqnarray}
It turns out that (\ref{E15a}) becomes
\begin{eqnarray}\label{E17}
\nonumber \Gamma(s+1)&=&    \Big({s\over e}\Big)^s\sqrt {2\pi s}\Big[ 1+{1\over 12s}+{1\over 288s^2}-{139\over 51840s^3}\\
&&  -{571\over 2488320s^4}+{163879\over 209018880s^5}  +o\big({1\over s^5}\big)\Big],\qquad \mbox{as}\, s\to +\infty.
\end{eqnarray}

\paragraph{Acknowledgments.}  This work was supported in part by NSFC (No.
61074047).

\bigskip

\noindent\textit{School of Mathematical Sciences, Fudan
University, Shanghai 200433, China \\
hwlou@fudan.edu.cn}

\newpage

\begin{center}\LARGE
Appendix: Methods to Calculate $a_n$.
\end{center}

First, we have $a_0={\sqrt 2\over 2}$.
\bigskip

\noindent\textbf{Method I.}
\begin{eqnarray*}
a_n&=& \lim_{v\to 0}{\displaystyle y(v)-\sum^{n-1}_{k=0}a_kv^k\over v^n}=\lim_{v\to 0^+}{\displaystyle y(v)-\sum^{n-1}_{k=0}a_kv^k\over v^n}\\
&=& \lim_{v\to 0^+}{\displaystyle {v\over f(v)}-\sum^{n-1}_{k=0}a_kv^k\over v^n}\\
&=& \lim_{x\to 0^+}{\displaystyle {\sqrt {x-\ln(1+x)}\over x}-\sum^{n-1}_{k=0}a_k (x-\ln(1+x))^{k\over 2}\over (x-\ln(1+x))^{n\over 2}}\\
&=& \lim_{x\to 0^+}{\displaystyle 2^{n\over 2}\Big({\sqrt {x-\ln(1+x)}\over x}-\sum^{n-1}_{k=0}a_k (x-\ln(1+x))^{k\over 2}\Big)\over x^n}, \qq n\geq 1.
\end{eqnarray*}

\noindent\textbf{Method II.}
It follows from
$$
v^2={v\over y(v)}-\ln \big(1+{v\over y(v)}\big)
$$
that
$$
{v\over v+y}\Big({v\over y}\Big)^\prime=2v.
$$
Then
$$
vy^\prime=y-2y^2-2 v y^2.
$$
For two smooth (not necessary analytic) functions $\vf(t)$ and $\psi(t)$, the MacLaurin's series of $\vf(t)\psi(t)$ is the Cauchy product of the
 MacLaurin's series of $\vf(t)$ and that of $\psi(t)$. Then we can establish that
\begin{eqnarray*}
a_n&=&\mbox{the coefficient of $v^n$ in the expansion of }\,\\
&& -{2\over n+2}\Big[\Big(\sum^{n-1}_{k=0}a_kv^k\Big)^3+v\Big(\sum^{n-1}_{k=0}a_kv^k\Big)^2\Big]\\
&=&-{2\over n+2}\Big[\sum_{k+j+l=n\atop 0\leq k,j,l<n}a_ka_ja_l+\sum_{k+j=n-1}a_ka_j\Big], \qq n\geq 1.
\end{eqnarray*}
More precisely,
for $n\geq 1$,  we have
\begin{eqnarray*}
&& a_0-\sum^n_{k=1}(k-1)a_kv^k\\
&=& 2\Big(\sum^n_{k=0}a_kv^k\Big)^3+2v\Big(\sum^n_{k=0}a_kv^k\Big)^2+o(v^n)\\
&=& 6a_0^2a_nv^n+2\Big(\sum^{n-1}_{k=0}a_kv^k\Big)^3+2v\Big(\sum^{n-1}_{k=0}a_kv^k\Big)^2+o(v^n)\\
&=& 3a_nv^n+2\Big(\sum^{n-1}_{k=0}a_kv^k\Big)^3+2v\Big(\sum^{n-1}_{k=0}a_kv^k\Big)^2+o(v^n), \qquad (v\to 0).
\end{eqnarray*}
That is,
$$
 a_0-\sum^n_{k=1}(k+2)a_kv^k=2\Big(\sum^{n-1}_{k=0}a_kv^k\Big)^3+2v\Big(\sum^{n-1}_{k=0}a_kv^k\Big)^2+o(v^n), \qquad (v\to 0).
$$
Therefore,
\begin{eqnarray*}
a_n&=&\mbox{the coefficient of $v^n$ in the expansion of }\,\\
&& -{2\over n+2}\Big[\Big(\sum^{n-1}_{k=0}a_kv^k\Big)^3+v\Big(\sum^{n-1}_{k=0}a_kv^k\Big)^2\Big].
\end{eqnarray*}

\bigskip

The advantage of Method I is that the convergence of the limit in calculating $a_{n+1}$ depends on the correctness of $a_0,a_1,\ldots,a_n$. Thus, we can
know whether $a_n$ is correct or not when we calculate $a_{n+1}$.

The advantage of Method II is that it brings possibility to calculate $a_n$ without the aid of computer.

We have:
$$
\begin{array}{lll}
&\displaystyle a_0={1\over 2}\sqrt 2, & \displaystyle a_1=-{1\over 3},\vspace{4mm}\\
& \displaystyle a_2={1\over 12}\sqrt 2,& \displaystyle a_3=-{4\over 135},\vspace{4mm}\\
&\displaystyle a_4={1\over 432}\sqrt 2, &\displaystyle a_5={4\over 2835},\vspace{4mm}\\
&\displaystyle a_6=-{139\over 194400}\sqrt 2, &\displaystyle  a_7={8\over 25515},\vspace{4mm}\\
&\displaystyle a_8=-{571\over 32659200}\sqrt 2,& \displaystyle a_9=-{1124\over 37889775},\vspace{4mm}\\
& \displaystyle a_{10}={163879\over 12345177600}\sqrt 2,\qq &\displaystyle a_{11}=-{41768\over 7388506125},\vspace{4mm}\\
& \displaystyle a_{12}={5246819\over 24443451648000}\sqrt 2,\qq &\displaystyle  a_{13}={43672\over 66496555125},\vspace{4mm}\\
& \displaystyle a_{14}=-{534703531\over 1906589228544000}\sqrt 2,\qq &\displaystyle  a_{15}={1459313264\over 12463116844303125},\vspace{4mm}\\
& \displaystyle a_{16}=-{4483131259\over 1372744244551680000}\sqrt 2,\qq &\displaystyle  a_{17}=-{10603947212\over 710397660125278125},\vspace{4mm}\\
& \displaystyle a_{18}={432261 921612371\over 69309856907414323200000}\sqrt 2,\qq &\displaystyle  a_{19}=-{49374413464\over 19180736823382509375},\vspace{4mm}\\
& \displaystyle a_{20}={6232523202521089\over 110618531624233259827200000}\sqrt 2.\qq &
\end{array}
$$

\bigskip

\begin{eqnarray*}
\nonumber \Gamma(s+1)&=&    \Big({s\over e}\Big)^s\sqrt {2\pi s}\Big[ 1+{1\over 12s}+{1\over 288s^2}-{139\over 51840s^3}-{571\over 2488320s^4}\\
&&\qq \qq\qq+{163879\over 209018880s^5}  +{5246819\over 75246796800s^6}- {534703531\over 902961561600s^7}\\
&& \qq \qq\qq-{4483131259\over 86684 309913600s^8}+{432261921612371\over 514904800886784000s^9}\\
&& \qq \qq\qq+{6232523202521089\over 86504006548979712000s^{10}}+o\big({1\over s^{10}}\big)\Big],\qquad \mbox{as}\, s\to +\infty.
\end{eqnarray*}

\end{document}